\documentclass[12pt]{article}
\usepackage{graphicx}        

\usepackage{amssymb, amsmath}
\usepackage{color}

\usepackage{textcomp}
\usepackage{psfrag}
\usepackage{epsfig}

\usepackage{latexsym}
\usepackage{amscd}
\newtheorem{theorem}{Theorem}
\newtheorem{proposition}{Proposition}
\newtheorem{lemma}{Lemma}

\newtheorem{remark}{Remark}

\def\P{{\cal P}}

\def\R{{\mathbb R}}

\def\ex{{\textrm{e}}}
\def\exx{\rm{\rm{e}}}
\def\po{\rm{\rm{Po}}}
\newcommand{\vertk}{\stackrel{\mbox{\scriptsize ${\cal D}$}}{\longrightarrow}}
\DeclareRobustCommand{\stirling2}{\genfrac\{\}{0pt}{}}

\def\PROB{{\mathbb P}}
\def\EXP{{\mathbb E}}

\def\IND{{\mathbb I}}

\def\qed{\hfill {\ \vrule width 1.5mm height 1.5mm \smallskip}}

\def\be{\begin{equation}}
\def\ee{\end{equation}}

\def\proof{\medskip  \par \noindent {\bf Proof.} \ \ }
\def\qedskip{\smallskip\noindent}
\def\qed{\hfill $\Box$ \qedskip}              

\begin{document}

\begin{titlepage}
\thispagestyle{empty}
\setcounter{page}{0}

\title{\textcolor{black}{The} limit distribution of the  maximum probability nearest neighbor ball}
\author{ L\'aszl\'o Gy\"orfi\thanks{Budapest University of Technology and Economics, gyorfi@cs.bme.hu}
\and Norbert Henze\thanks{Karlsruhe Institute of Technology, henze@kit.edu}
\and Harro Walk\thanks{Universit\"at Stuttgart, walk@mathematik.uni-stuttgart.de}}

\maketitle

\begin{abstract}
Let $X_1,\ldots,X_n$ be independent random points drawn from an absolutely continuous
probability measure with density $f$ in $\R^d$. Under mild conditions on $f$, we derive a Poisson limit
theorem for the number of large probability nearest neighbor balls. Denoting by $P_n$ the maximum probability measure
of nearest neighbor balls, this limit theorem implies a Gumbel extreme value distribution for $nP_n - \ln n$ as $n \to \infty$.
Moreover, we derive a tight upper bound on the upper tail of the distribution of $nP_n - \ln n$, which does not depend on $f$.
\end{abstract}

\bigskip


\bigskip


\end{titlepage}

\section{\textcolor{black}{Introduction} }

Let $X,X_1,\ldots,X_n,\ldots$ be independent, identically distributed (i.i.d.) random
vectors  taking values in $\R^d$. We assume throughout the paper that the distribution of
$X$, which is denoted by $\mu$, has a density $f$ with respect to Lebesgue measure $\lambda$.

Writing $\|\cdot \|$ for the Euclidean norm on $\R^d$, put
\[
R_{i,n}:=\min_{j\ne i, j\le n} \|X_i-X_j\|,
\]
and let
\[
P_n:=\max_{1\le i \le n}\mu \{S(X_i,R_{i,n})\}
\]
denote the maximum probability of the nearest neighbor (NN) balls,
where $S(x,r) := \{y \in \R^d: \|y-x\| \le r\}$ stands for the closed ball with center $x$ and radius $r$.
This paper deals with both the finite-sample and the asymptotic distribution of
\[
nP_n - \ln n,
\]
as $n \to \infty$.

There is a huge related literature for Poisson sample size.
Let $N$ be a random variable that is independent of $X_1,X_2,\dots$ and has  a Poisson distribution with $\mathbb{E}(N) = n$.
Then
\begin{align}
\label{poi}
X_1,\dots, X_N
\end{align}
is a non-homogeneous Poisson process with intensity function $nf$.
For the nucleuses  $X_1,\dots, X_N$,
$\widetilde A_n(X_j)$ denotes the Voronoi cell around $X_j$, and $\widehat r_j$ and $\widehat R_j$ stand for the inscribed and circumscribed radii of $\widetilde A_n(X_j)$, respectively,
i.e., we have
\[
\widehat r_j=\sup\{r>0:S(X_j,r)\subset \widetilde A_n(X_j)\}
\]
and
\[
\widehat R_j=\inf\{r>0: \widetilde A_n(X_j)\subset S(X_j,r)\}.
\]

If  $X_1,X_2,\dots$ are i.i.d. uniformly distributed on the unit cube $[0,1]^{d}$, then
(2a) and (2c) of Theorem 1 in Calka and  Chenavier \cite{CaCh14} read
\begin{align*}
\lim_{n\to \infty}\PROB\left( 2^{d}n\lambda \left\{S\left(0,\max_{1\le j\le N}\widehat r_j\right)\right\} -\ln n\le y\right)=G(y)
\end{align*}
and
\begin{align*}
\lim_{n\to \infty}\PROB\left( n\lambda \left\{ S\left( 0,\max_{1\le j\le N}\widehat R_j)\right)\right\} -\ln \left(\alpha_d n\left(\ln n\right)^{d-1}\right)
\le y\right)=G(y),
\end{align*}
$y\in \R$. Here, $\alpha_d>0$ is a universal constant, and
\[
G(y)=\exp(-\exp(-y))
\]
denotes the distribution function of the Gumbel extreme value distribution.

The paper is organized as follows. In Section \ref{sec2} we study the distribution of $nP_n -\ln n$.
Theorem \ref{cor} is on a universal and tight bound on the upper tail of $nP_n -\ln n$. Under mild conditions on the density,
 Theorem \ref{G2d} shows that the number of exceedances of nearest neighbor ball probabilities over a certain sequence of thresholds
 has an asymptotic Poisson distribution as $n \to \infty$. As a consequence, the limit distribution of $nP_n -\ln n$
 is the Gumbel extreme value distribution. Theorem \ref{Poisson} in Section \ref{sec3} is the extension of Theorem \ref{cor} for Poisson sample size.
 All proofs are presented in Section \ref{sec4}. The main tool for proving Theorem \ref{G2d} is a novel Poisson limit theorem for sums of indicators of exchangeable events,
 which is formulated as Proposition \ref{henze}.  The final section sheds some light on a technical condition on $f$ that is used in the proof of the main result.

 Although there is a weak dependence between the probabilities of nearest neighbor balls, a main message of this paper is that one
 can neglect this dependence when looking for the limit distribution of the maximum probability.

\section{The maximum nearest neighbor ball }
\label{sec2}

Under the assumption that the density $f$ is sufficiently smooth and bounded away from zero,
Henze \cite{Hen81} and  \cite{Hen82} derived the limit distribution of the maximum {\em approximate} probability measure
\begin{align}
\label{Appr}
\max_{1\le i\le n}f(X_i)R_{i,n}^dv_d
\end{align}
of NN-balls. Here, $v_d= \pi^{d/2}/\Gamma(1+d/2)$ stands for the volume of the unit ball in $\R^d$.

In the following, we consider the number of points among $X_1,\ldots,X_n$ for which the probability content of the nearest
neighbor ball exceeds some (large) threshold. To be more specific, we fix $y \in \R$ and consider the random variable
\[
C_n := \sum_{i=1}^n \IND\big{\{}n \mu\{S(X_i,R_{i,n})\} > y + \ln n\big{\}},
\]
where $\IND\{ \cdot \}$ denotes the indicator function.  Writing "$\vertk$" for convergence in distribution,
we will show that, under some conditions on the density $f$,
\[
C_n \vertk Z \ \ \textrm{ as } n \to \infty,
\]
where $Z$ is a random variable with the Poisson distribution Po$(\exp(-y))$. Now, $C_n =0$ if, and only if, $n P_n - \ln n \le y$,
and it follows that
\begin{align}
\label{appNN}
\lim_{n \to \infty} \PROB\left(nP_n - \ln n \le y\right) = \PROB(Z=0) = G(y),\quad y\in \R.
\end{align}

\bigskip
Since $1-G(y) \le \exp(-y)$ if $y \ge 0$, (\ref{appNN}) implies
\begin{align}
\label{uNN}
\limsup_{n\to \infty}\PROB\left(nP_n - \ln n \ge y\right)\le  \ex^{-y},\quad y\ge 0.
\end{align}
Our first result is a non-asymptotic upper bound on the upper tail of the distribution of $nP_n - \ln n$. This bound holds without any condition on the density and thus entails  (\ref{uNN}) universally.
\begin{theorem}
\label{cor}
Without any restriction on the density $f$, we have
\begin{align}
\label{cNN}
\PROB\left(nP_n - \ln n \ge y\right)\le \exp\left(-\frac{n-1}{n}y+\frac{\ln n}{n}\right)\IND\{y\le n-\ln n\},\quad y\in \R.
\end{align}
\end{theorem}

Theorem \ref{cor} implies a non-asymptotic upper bound on the mean of $nP_n - \ln n $, since
\begin{align*}
\EXP\left[nP_n - \ln n \right]
&\le
\EXP\left[(nP_n - \ln n )^+\right]\\
&=
\int_0^{\infty}
\PROB\left(nP_n - \ln n \ge y\right)\, \textrm{d}y\\
&\le
\int_0^{\infty}\exp\left(-\frac{n-1}{n}y+\frac{\ln n}{n}\right)\, \textrm{d}y\\
&=
\frac{n}{n-1}\exp\left(\frac{\ln n}{n}\right).
\end{align*}
Notice that this upper bound approaches $1$ for large $n$, and  that the mean of the standard Gumbel distribution is the Euler-Mascheroni
constant, which is $-\int_0^{\infty}\ex^{-y}\ln y\, \textrm{d}y=0.5772\dots$

\bigskip

Recall that the support of $\mu$ is defined by
\[
\textrm{supp}(\mu) := \{x \in \R^d: \mu\{S(x,r)\} > 0 \textrm{ for each } r >0\},
\]
i.e., the support of $\mu$ is the smallest closed set in $\R^d$ having $\mu$-measure one.

\begin{theorem}
\label{G2d}
Assume there are $\beta \in (0,1)$, $c_{max}<\infty$ and $\delta >0$ such that, for any $r,s >0$ and any $x,z \in$ supp$(\mu)$ with
$\|x-z\| \ge \max\{r,s\}$ and  $\mu\left( S(x,r)\right) = \mu\left( S(z,s)\right)\le \delta$, one has
\begin{align}
\label{INT}
\frac{\mu\left( S(x,r)\cap   S(z,s)\right)}{\mu\left( S(z,s)\right)}
\le \beta
\end{align}
and
\begin{align}
\label{L}
\mu\left( S(z, 2s)\right)\le c_{max}\mu\left( S(z,s)\right).
\end{align}
Then
\begin{align}\label{poissconv}
\sum_{i=1}^n \IND\big{\{}n \mu\{S(X_i,R_{i,n})\} > y + \ln n\big{\}} \vertk {\po}(\exp(-y)), \quad y \in \R,
\end{align}
and hence
\begin{align}
\label{NN}
\lim_{n\to \infty}\PROB\left( nP_n - \ln n \le y  \right)=G(y),\quad y\in \R.
\end{align}
\end{theorem}

\begin{remark} It is easy to see that
\label{Corr}
(\ref{INT}) and (\ref{L}) hold if the density is both bounded from above by $f_{max}$ and bounded away from zero by $f_{min}>0$. Indeed, putting
\[
\beta := 1 - \frac{1}{2} \cdot \frac{f_{min}}{f_{max}}, \quad c_{max} := 2^d \cdot \frac{f_{max}}{f_{min}},
\]
we have
\begin{align*}
\frac{\mu\left( S(x,r)\cap   S(z,s)\right)}{\mu\left( S(z,s)\right)}
&=
1-\frac{\mu\left( S(z,s)\setminus S(x,r)\right)}{\mu\left( S(z,s)\right)}\\
&\le
1-\frac{f_{min}\, \lambda\left( S(z,s)\setminus S(x,r)\right)}{f_{max}\, \lambda\left( S(z,s)\right)}\\
&\le
\beta
\end{align*}
and
\begin{align*}
\mu\left( S(z, 2s)\right)
&\le
f_{max} \, \lambda\left( S(z, 2s)\right)\\
&=
f_{max} \, 2^d\lambda\left( S(z, s)\right)\\
&\le
c_{max} \, \mu\left( S(z, s)\right).
\end{align*}

A challenging problem left is to weaken the conditions of Theorem \ref{G2d} or to prove that (\ref{poissconv}) and (\ref{NN}) hold
without any conditions on the density. We believe that such universal limit results are possible, because the summands in (\ref{poissconv}) are identically distributed, and their
distribution does not depend on the actual density. More discussion on condition (\ref{INT}) is given in Section
\ref{secdiscuss}.
\end{remark}

\section{The maximum nearest neighbor ball for a non-homogeneous Poisson process}
\label{sec3}

In this section we consider the non-homogeneous Poisson process $X_1,\dots, X_N$  defined by (\ref{poi}).
Putting
\[
\widetilde R_{i,n}:=\min_{j\ne i, j\le N} \|X_i-X_j\|
\]
and
\[
\widetilde P_n=\max_{1\le i \le N}\mu \{S(X_i,\widetilde R_{i,n}) \},
\]
the following result is the Poisson-analogue to Theorem \ref{cor}.
\begin{theorem}
\label{Poisson}
Without any restriction on the density $f$ we have
\begin{align*}
\PROB\left(n\widetilde P_n - \ln n \ge y\right)
&\le  {{\exx}^{-y}} \exp\left(\frac{(y+ \ln n)^2}{n}\right),\quad y\in \R.
\end{align*}
\end{theorem}

\section{Proofs }
\label{sec4}

\bigskip
\noindent
{\bf Proof of Theorem \ref{cor}.}
Since the right hand side of (\ref{cNN}) is larger than $1$ if $y <0$, we take $y \ge 0$ in what  follows.
Moreover, in view of $P_n\le 1$ the left hand side of (\ref{cNN}) vanishes if $y > n-\ln n$. We therefore assume without loss of generality
that
\begin{align}
\label{z}
\frac{y+ \ln n}{n}\le 1.
\end{align}
For a fixed $x \in \R^d$, let
\begin{equation}\label{defhx}
H_x(r):=\PROB \left(\|x-X\|\le r\right), \quad r \ge 0,
\end{equation}
\textcolor{black}{be} the distribution function of $\|x-X\|$. By the probability integral transform (cf. Biau and Devroye \cite{BiDe15}, p. 8),
the random variable
\[
H_x(\|x-X\|)=\mu \{ S(x,\|x-X\|) \}
\]
is uniformly distributed on $[0,1]$. We thus have
\begin{align}
\label{unif}
\mu \left\{ S(x,H^{-1}_x(p)) \right\}= p, \quad 0 < p < 1,
\end{align}
where $H^{-1}_x(p) = \inf\{r: H_x(r) \ge p\}$. It follows that
\begin{align*}
\PROB\left(nP_n - \ln n \ge y\right)
&=
\PROB\left(n\max_{1\le i \le n}\mu \{ S(X_i,R_{i,n})\} - \ln n \ge y\right)\\
&\le
n\, \PROB\left(n\mu \{ S(X_1,R_{1,n}) \} - \ln n \ge y\right)\\
&=
n\, \PROB\left(\mu \{ S(X_1,R_{1,n}) \}  \ge \frac{y+ \ln n}{n}\right)\\
&=
n\, \PROB\left(\min_{2\le j\le n}\mu \{ S(X_1,\|X_1-X_j\|) \}  \ge \frac{y+ \ln n}{n}\right).
\end{align*}
Now, (\ref{z}) implies
\begin{align*}
\PROB\left(nP_n - \ln n \ge y\right)
&\le
n\, \EXP \left[\PROB \left( \min_{2\le j\le n}\mu \{S(X_1,\|X_1-X_j\|) \}  \ge \frac{y+ \ln n}{n} \Big{|} X_1\right)\right]\\
&=
n\left(1- \frac{y+ \ln n}{n}\right)^{n-1}\\\
&\le
n \exp\left(- \frac{(y+ \ln n)(n-1)}{n}\right)\\
&=
\exp\left(-\frac{n-1}{n}y+\frac{\ln n}{n}\right).
\end{align*}
\qed

\bigskip
\noindent
{\bf Proof of Theorem \ref{Poisson}.}
We again assume (\ref{z}) in what follows. By conditioning on $N$, we have
\begin{align*}
\PROB\left(n\widetilde P_n - \ln n \ge y\right)
&=
\sum_{k=1}^{\infty}\PROB\left(n\widetilde P_n - \ln n \ge y\mid N=k\right)\PROB\left(N=k\right)\\
&=
\sum_{k=1}^{\infty}\PROB\left(nP_k - \ln n \ge y\right)\PROB\left(N=k\right).
\end{align*}
Putting $y_n:= (y+ \ln n)/n$, we obtain
\begin{align*}
\PROB\left(nP_k - \ln n \ge y\right)
&=
\PROB\left(kP_k - \ln k \ge k y_n-\ln k\right),
\end{align*}
and Theorem \ref{cor} implies
\begin{align*}
\PROB \! \left(kP_k - \ln k \ge k y_n -\ln k\right)
&\le
\exp \! \left(-\frac{k\! -\! 1}{k}\left( k y_n -\ln k \right)+\frac{\ln k}{k}\right)\\
&=
\exp \! \left(-(k-1)y_n +\ln k\right).
\end{align*}
It follows that
\begin{align*}
\PROB\left(n\widetilde P_n - \ln n \ge y\right)
&\le
\sum_{k=1}^{\infty} \exp\left(-(k-1)y_n+\ln k\right) \PROB\left(N=k\right)\\
&=
\ex^{y_n-n} \sum_{k=1}^{\infty}k\left(\ex^{-y_n}\right)^k\frac{n^k}{k!}\\
&=
\ex^{y_n-n} \sum_{k=1}^{\infty}\frac{\left(n\ex^{-y_n}\right)^k}{(k-1)!}\\
&=
n\ex^{y_n-n-y_n} \exp\left(n\ex^{-y_n}\right)\\
&=
n \exp(-n(1-\ex^{-y_n})).
\end{align*}
Since $z \ge 0$ entails $\ex^{-z}\le 1-z+z^2$, we finally obtain
\begin{align*}
\PROB\left(n\widetilde P_n - \ln n \ge y\right)
&\le
n \exp\left(-n(y_n-y_n^2)\right) =
\ex^{-y}\exp\left(\frac{(y+ \ln n)^2}{n}\right).
\end{align*}
\qed

\bigskip
\noindent
The main tool in the proof Theorem \ref{G2d} is the following result.
\smallskip
\begin{proposition}
\label{henze}
For each $n \ge 2$, let $A_{n,1}, \ldots, A_{n,n}$ be exchangeable events, and let
\[
Y_n := \sum_{j=1}^n \IND\{A_{n,j}\}.
\]
If, for some $\nu \in (0,\infty)$,
\begin{align}
\label{cond}
\lim_{n\to \infty}n^k \PROB\left(A_{n,1} \cap \ldots \cap A_{n,k}\right) = \nu^k \quad \textrm{for each } \ k \ge 1,
 \end{align}
then
\[
Y_n \vertk Y \quad \textrm{ as } \ n \to \infty,
\]
where $Y$ has the Poisson distribution \po$(\nu)$.
\end{proposition}

\noindent
{\bf Proof.} The proof uses the method of moments, see, e.g., \cite{Billi}, Section 30.
Putting
\[
S_{n,k} = \sum_{1 \le i_1 < \cdots < i_k \le n} \PROB\left(A_{n,i_1} \cap \ldots \cap A_{n,i_k}\right), \quad k \in \{1,\ldots,n\},
\]
and writing $Z^{(k)} = Z(Z-1) \cdots (Z-k+1)$ for the $k$th descending factorial of a random variable $Z$, we have
\[
\mathbb{E}\left[Y_n^{(k)}\right] = k! \, S_{n,k}.
\]
Since $A_{n,1}, \ldots, A_{n,n}$ are exchangeable, (\ref{cond}) implies
\[
\lim_{n\to \infty} \mathbb{E}\left[Y_n^{(k)}\right] = \nu^k, \quad k \ge 1.
\]
Now, $\nu^k = \mathbb{E}\left[Y^{(k)}\right]$, where $Y$ has the Poisson distribution Po$(\nu)$. We thus have
\begin{equation}\label{momconv}
\lim_{n \to \infty} \mathbb{E}\left[Y_n^{(k)}\right] =   \mathbb{E}\left[Y^{(k)}\right], \quad k \ge 1.
\end{equation}
Since
\[
Y_n^k = \sum_{j=0}^k \stirling2{k}{j} Y_n^{(j)},
\]
where $\stirling2{k}{0}, \ldots , \stirling2{k}{k}$ denote Stirling numbers of the second kind (see, e.g., \cite{GKP}, p. 262),
(\ref{momconv}) entails $\lim_{n \to \infty} \mathbb{E}[Y_n^k] = \mathbb{E}[Y^k]$ for each $k \ge 1$. Since the distribution
of $Y$ is uniquely determined by the sequence of moments $(\mathbb{E}[Y^k])$, $k \ge 1$, the assertion follows.
\qed

\vspace*{5mm}

\noindent {\bf Proof of Theorem \ref{G2d}.}
Fix $y \in \R$. In what follows, we will verify (\ref{cond})
for
\[
A_{n,i} := \big{\{} n\mu \{S(X_i,R_{i,n} )\} \ge y + \ln n\big{\}}, \quad i \in \{1,\ldots,n\},
\]
and $\nu = \exp(-y)$. Throughout the proof we tacitly assume
\[
0 < y_n := \frac{y+\ln n}{n} < 1.
\] 
This assumption entails no loss of generality since $n$ tends to infinity.
With $H_x(\cdot)$ given in (\ref{defhx}), we put
\[
R^*_{i,n}:=H_{X_i}^{-1}((y+\ln n)/n), \quad i \in \{1,\ldots,n\}.
\]
For the special case $k=1$, conditioning on $X_1$ and (\ref{unif}) yield
\begin{align*}
n\,  \PROB \left(A_{n,1}\right)
& =  n\,  \PROB\left(\mu\left(S(X_1,R_{1,n}) \right) \ge y_n \right)\\
& = n\,  \mathbb{E} \Big{[} \PROB\left(\mu\left(S(X_1,R_{1,n}) \right) \ge y_n  \mid X_1 \right) \Big{]}\\
& = n\, \mathbb{E} \Big{[} \left(1- \mu\left(S(X_1, H_{X_1}^{-1}(y_n))\right) \right)^{n-1}  \Big{]}\\
& = n \, \left(1 - \frac{y + \ln n}{n}\right)^{n-1}.
\end{align*}
Using the inequalities $1-1/t \le \ln t \le t-1$ gives $\lim_{n \to \infty} n\, \PROB(A_{n,1}) = \ex^{-y}$.
Thus (\ref{cond}) is proved for $k=1$, remarkably without any condition on the underlying density $f$.
We now assume $k \ge 2$ and put
\[
\widetilde R_{i,k,n}:=\min_{k+1\le j\le n} \|X_i-X_j\|, \quad
r_{i,k}:=\min_{j\ne i, j\le k} \|X_i-X_j\|.
\]
Then
\[
R_{i,n}=\min\{\widetilde R_{i,k,n},r_{i,k}\},
\]
and because of $\widetilde R_{i,k,n}\to 0$ $\mathbb{P}$-almost surely as $n \to \infty$,
it follows that, on a set of probability 1,
\[
R_{i,n}=\widetilde R_{i,k,n} \quad \textrm{for each } i \in \{1,\ldots,k\}
\]
if $n$ is large enough.
Conditioning on $X_1,\ldots,X_k$  we  have
\begin{align*}
&\PROB\left(\cap_{i=1}^k A_{n,i}\right)\\
&= \PROB\left(\cap_{i=1}^k \{ \mu (S(X_i,\min\{\widetilde R_{i,k,n},r_{i,k}\}))\ge y_n\}\right)\\
&= \PROB\left(\cap_{i=1}^k \{ \mu (S(X_i,\widetilde R_{i,k,n} ))\ge y_n, \mu (S(X_i,r_{i,k} ))\ge y_n\}\right)\\
&= \EXP\bigg{[}\PROB\left(\cap_{i=1}^k \{ \mu (S(X_i,\widetilde R_{i,k,n}))\ge y_n, \mu (S(X_i,r_{i,k} ))\ge y_n\}
\mid X_1, \dots ,X_k\right)\bigg{]}\\
&= \EXP\bigg{[}\PROB \! \left( \! \cap_{i=1}^k \{ \mu (S(X_i,\widetilde R_{i,k,n} )) \! \ge \! y_n\} \! \mid \! X_1, \dots ,X_k\!  \right)
\prod_{i=1}^k \IND\{\mu (S(X_i,r_{i,k} )) \! \ge \! y_n\}\! \bigg{]}\! .
\end{align*}
Furthermore, we obtain
\begin{align*}
&\PROB\left(\cap_{i=1}^k \{ \mu (S(X_i,\widetilde R_{i,k,n} ))\ge y_n\}\mid X_1, \dots ,X_k\right)\\
&=  \PROB\left(\cap_{i=1}^k \{ \widetilde R_{i,k,n}\ge H_{X_i}^{-1}(y_n)\}\mid X_1, \dots ,X_k\right)\\
&=  \PROB\left(\cap_{i=1}^k \{ \widetilde R_{i,k,n}\ge R^*_{i,n}\}\mid X_1, \dots ,X_k\right)\\
&=  \PROB\Big{(}X_{k+1},\dots ,X_n \notin \cup_{i=1}^k S(X_i, R^*_{i,n})\mid X_1, \dots ,X_k\Big{)}\\
&=  \left(1-\mu\left( \cup_{i=1}^k S(X_i, R^*_{i,n})\right)\right)^{n-k}.
\end{align*}
Notice that we have the obvious lower bound
\begin{align*}
& n^k \left(1-\mu\left( \cup_{i=1}^k S(X_i, R^*_{i,n})\right)\right)^{n-k}
\prod_{i=1}^k \IND\{\mu (S(X_i,r_{i,k} ))\ge y_n\}\\
&\ge n^k \left(1-\sum_{i=1}^k\mu\left(  S(X_i, R^*_{i,n})\right)\right)^{n-k}
\prod_{i=1}^k \IND\{\mu (S(X_i,r_{i,k} )\ge y_n\}\\
&= n^k \left(1- k\, \frac{y+\ln n}{n} \right)^{n-k}\prod_{i=1}^k \IND\{\mu (S(X_i,r_{i,k} ))\ge y_n\}.
\end{align*}
Since the latter converges almost surely to $\ex^{-ky}$ as $n \to \infty$,
Fatou's lemma implies
\begin{align*}
& \liminf_{n\to \infty} n^k \PROB\left(\cap_{i=1}^k A_{n,i}\right)\\
&\textcolor{blue}{=} \liminf_{n\to \infty} \EXP\Big{[} n^k \left(1-\mu\left( \cup_{i=1}^k S_{X_i, R^*_{i,n}}\right)\right)^{n-k}
\prod_{i=1}^k \IND\{\mu (S(X_i,r_{i,k} ))\ge y_n\}\Big{]}\\
&\ge \EXP\Big{[}\liminf_{n\to \infty} n^k \left(1-\mu\left( \cup_{i=1}^k S_{X_i, R^*_{i,n}}\right)\right)^{n-k}
\prod_{i=1}^k \IND\{\mu (S(X_i,r_{i,k} ))\ge y_n\}\Big{]}\\
&=\ex^{-ky}.
\end{align*}
It thus remains to show
\begin{equation}\label{upestim}
\limsup_{n\to \infty} n^k \PROB\left(\cap_{i=1}^k A_{n,i}\right) \le \ex^{-ky}.
\end{equation}
Let $D_n$ be the event that the balls $S(X_i, R^*_{i,n})$, $i=1,\dots ,k$, are pairwise
disjoint. Putting
\[
\IND_{n,k} := \prod_{i=1}^k \IND\{\mu (S(X_i,r_{i,k} ))\ge y_n\},
\]
we have
\begin{align*}
& \limsup_{n\to \infty} n^k \PROB(A_{n,1}\cap \ldots \cap A_{n,k})\\
& = \limsup_{n \to \infty} n^k \, \EXP\bigg{[}\left(1-\mu\left( \cup_{i=1}^k S(X_i, R^*_{i,n})\right)\right)^{n-k}
 \IND_{n,k}\bigg{]}\\
&\le\limsup_{n \to \infty} n^k \, \EXP\bigg{[}\exp\left(-(n-k)\mu\left( \cup_{i=1}^k S(X_i, R^*_{i,n})\right)\right)
 \IND_{n,k}\bigg{]}\\
&\le
\limsup_{n \to \infty}  n^k \, \EXP\bigg{[}\exp\left(-(n-k)k\frac{y+\ln n}{n} \right)
 \IND\{D_n\}\bigg{]}\\
&+\limsup_{n \to \infty} n^k \, \EXP\bigg{[}\exp\left(-(n-k)\mu\left( \cup_{i=1}^k S(X_i, R^*_{i,n})\right)\right)
\IND\{D^c_n\} \IND_{n,k}\bigg{]}\\
&\le
\ex^{-ky}
+\limsup_{n \to \infty}  n^k \, \EXP\bigg{[}\exp\left(-(n-k)\mu\left( \cup_{i=1}^k S(X_i, R^*_{i,n})\right)\right)
\IND\{D^c_n\} \IND_{n,k}\bigg{]}.
\end{align*}
It thus remains to show
\begin{align}
\label{ndis}
\lim_{n \to \infty} n^k\,  \EXP\big{[}\exp\left(-(n-k)\mu\left( \cup_{i=1}^k S(X_i, R^*_{i,n})\right)\right)
\IND\{D^c_n\} \IND_{n,k} \big{]}
=0.
\end{align}
Under some additional smoothness conditions on the density, Henze \cite{Hen81} verified (\ref{ndis}) for the related problem
of finding the limit distribution of the random variable figuring in (\ref{Appr}). By analogy with his way of proof, we
introduce an equivalence relation on the set $\{1,\dots ,k\}$ as follows:  An equivalence class consists of a singleton $\{i\}$ if
\[
S(X_i, R^*_{i,n})\cap S(X_j, R^*_{j,n})=\emptyset
\]
for each $j\ne i$. Otherwise, $i$ and $j$ are called equivalent if there is a subset $\{i_1,\dots ,i_\ell\}$ of $\{1,\dots ,k\}$ such that $i=i_1$ and $j=i_\ell$ and
\[
S(X_{i_m}, R^*_{i_m,n})\cap S(X_{i_{m+1}}, R^*_{i_{m+1},n})\ne\emptyset
\]
for each $m\in \{1,\dots ,\ell-1\}$.
Let $\P=\{Q_1,\dots ,Q_q\}$
be a partition of $\{1,\dots,k\}$, and denote by $E_u$ the event that $Q_u$ forms an equivalence class.
For the event $D_n$, the partition
$
\P_0:=\{\{1\},\dots ,\{k\}\}
$
is the trivial one, while
on the complement $D^c_n$ any partition $\P$ is non-trivial, which means that $q<k$.
In order to prove (\ref{ndis}), we have to show that
\begin{align}
\label{Q}
\limsup_{n\to \infty}  n^k \EXP\bigg{[}\exp\left(-(n-k)\mu\left( \cup_{i=1}^k S(X_i, R^*_{i,n})\right)\right)
\IND_{n,k}
\prod_{u=1}^q\IND\{E_u\} \bigg{]}
&=0
\end{align}
for each  non-trivial partition $\P$. Since balls that belong to different equivalence classes are disjoint,
we have
\begin{align*}
\mu\left( \cup_{i=1}^k S(X_i, R^*_{i,n})\right)\prod_{u=1}^q\IND\{E_u\}
&=
\mu\left( \cup_{u=1}^q \cup_{i\in Q_u}S(X_i, R^*_{i,n})\right)\prod_{u=1}^q\IND\{E_u\}\\
&=
\sum_{u=1}^q\mu\left(  \cup_{i\in Q_u}S(X_i, R^*_{i,n})\right)\prod_{u=1}^q\IND\{E_u\}.
\end{align*}
Writing $|B|$ for the number of elements of a finite set $B$, it follows that
\begin{align*}
& n^k \exp\left(-(n-k)\mu\left( \cup_{i=1}^k S(X_i, R^*_{i,n})\right)\right)
\prod_{i=1}^k \IND\{\mu (S(X_i,r_{i,k} ))\ge y_n\}
\prod_{u=1}^q\IND\{E_u\} \\
&\le
\ex^k
\prod_{u=1}^q n^{|Q_u|} \prod_{u=1}^q\ex^{-n\mu\left( \cup_{i\in Q_u} S(X_i, R^*_{i,n})\right)}
\prod_{u=1}^q\prod_{i\in Q_u}\! \IND\{\mu (S(X_i,r_{i,k} ))\! \ge \!  y_n\}
\prod_{u=1}^q\! \IND\{\! E_u\! \}\\
&=
 \ex^k
\prod_{u=1}^q \left(n^{|Q_u|} \ex^{-n\mu\left( \cup_{i\in Q_u} S(X_i, R^*_{i,n})\right)}
\prod_{i\in Q_u}\IND\{\mu (S(X_i,r_{i,k} ))\ge y_n\}
\IND\{E_u\}\right).
\end{align*}
In view of independence, we have
\begin{align*}
& n^k \EXP\bigg{[}\ex^{-n\mu\left( \cup_{i=1}^k S(X_i, R^*_{i,n})\right)}
\prod_{i=1}^k \IND\{\mu (S(X_i,r_{i,k} ))\ge y_n\}
\prod_{u=1}^q\IND\{E_u\} \bigg{]}\\
&=
\prod_{u=1}^q\EXP\bigg{[}n^{|Q_u|} \ex^{-n\mu\left( \cup_{i\in Q_u} S(X_i, R^*_{i,n})\right)}
\prod_{i\in Q_u}\IND\{\mu (S(X_i,r_{i,k} ))\ge y_n\}
\IND\{E_u\} \bigg{]}.
\end{align*}
Thus, (\ref{Q}) is proved if we can show
\begin{align*}
\lim_{n\to \infty} \EXP\bigg{[}n^{|Q_u|} \ex^{-n\mu\left( \cup_{i\in Q_u} S(X_i, R^*_{i,n})\right)}
\prod_{i\in Q_u}\IND\{\mu (S(X_i,r_{i,k} ))\ge y_n\}
\IND\{E_u\} \bigg{]}
& = 0
\end{align*}
for each $u$ with $2 \le |Q_u|< k$. Without loss of generality assume
\[
Q_u=\{1,\dots ,|Q_u|\}.
\]
Then
\begin{align*}
&\cap_{i=1}^{|Q_u|} \{\mu (S(X_i,r_{i,k} ))\ge y_n\}\\
&=\cap_{i=1}^{|Q_u|} \Big{\{}\mu (S(X_i,\min_{j\ne i, j\le |Q_u|} \|X_i-X_j\| ))\ge y_n\Big{\}}\\
&=\cap_{i=1}^{|Q_u|} \Big{\{}\min_{j\ne i, j\le |Q_u|}\mu (S(X_i, \|X_i-X_j\| ))\ge y_n\Big{\}}\\
&=\cap_{i=1}^{|Q_u|} \cap_{j\ne i, j\le |Q_u|}\{\mu (S(X_i, \|X_i-X_j\| ))\ge y_n\}\\
&=\cap_{i=1}^{|Q_u|} \cap_{j\ne i, j\le |Q_u|}\{\|X_i-X_j\|\ge H_{X_i}^{-1}(y_n)\}\\
&=\cap_{ i,j\le |Q_u|, i\ne j}\{\|X_i-X_j\|\ge \max(R^*_{i,n},R^*_{j,n})\},
\end{align*}
and we obtain
\begin{align*}
&n^{|Q_u|} \ex^{-n\mu\left( \cup_{i\in Q_u} S(X_i, R^*_{i,n})\right)}
\prod_{i\in Q_u}\IND\{\mu (S_{X_i,r_{i,k} })\ge y_n\}
\IND\{E_u\} \\
&=
n^{|Q_u|} \ex^{-n\mu\left( \cup_{i=1}^{|Q_u|} S(X_i, R^*_{i,n})\right)}
\IND\{\cap_{ i,j\le |Q_u|, i\ne j}\{\|X_i\! -\! X_j\| \! \ge \! \max\{R^*_{i,n},\! R^*_{j,n}\}\! \}\! \}
\IND\{\! E_u\! \} \\
&\le
n^{|Q_u|} \ex^{-n\mu\left( \cup_{i=1}^2 S(X_i, R^*_{i,n})\right)}
\IND\{\cap_{ i,j\le |Q_u|, i\ne j}\{\|X_i\! -\! X_j\| \! \ge \! \max\{R^*_{i,n},\! R^*_{j,n}\}\! \}\! \}
\IND\{\! E_u\! \}.
\end{align*}
Now, condition (\ref{INT}) implies
\begin{align*}
&n\mu\left( \cup_{i=1}^2 S(X_i, R^*_{i,n})\right)\\
&=
n\mu\left( S(X_1, R^*_{1,n})\right)
+
n\mu\left( S(X_2, R^*_{2,n})\right)
-
n\mu\left( S(X_2, R^*_{2,n})\cap  S(X_1, R^*_{1,n})\right)\\
&=
n\, \frac{y+\ln n}{n}\left(2-\frac{\mu\left( S(X_2, R^*_{2,n})\cap  S(X_1, R^*_{1,n})\right)}{\mu\left( S(X_2, R^*_{2,n})\right)}\right)\\
&\ge
(y+\ln n)(2-\beta)\\
&=:
(y+\ln n)(1+\varepsilon)
\end{align*}
(say). Notice that $\varepsilon >0$ since $0 < \beta < 1$.
Thus,
\begin{align*}
&n^{|Q_u|} \EXP\bigg{[}\ex^{-n\mu\left( \cup_{i\in Q_u} S(X_i, R^*_{i,n})\right)}
\prod_{i\in Q_u}\IND\{\mu (S(X_i,r_{i,k} ))\ge y_n\}
\IND\{E_u\}\bigg{]} \\
&\le
n^{|Q_u|} \ex^{-(y+\ln n)\left(1+\varepsilon \right)}
\EXP\bigg{[}\IND\{\cap_{ i,j\le |Q_u|, i\ne j}\{\|X_i\! -\! X_j\|\! \ge \! \max\{R^*_{i,n},R^*_{j,n}\}\! \}\! \}
\IND\{\! E_u\! \}\bigg{]}\\
&=
\textrm{O}\left(n^{|Q_u|-1-\varepsilon}\right)
\PROB\left(E_u\right).
\end{align*}
In order to bound $\PROB(E_u)$ we need the following lemma:
\begin{lemma}
\label{harro}
On $E_u$ there is a random integer $L\in\{1,\dots ,|Q_u|\}$ depending on $X_1,\dots ,X_{|Q_u|}$ such that $Q_u\setminus \{L\}$ forms an equivalence class.
\end{lemma}
\proof
Let $m:= |Q_u|$. Regard $X_1,\ldots,X_m$ as vertices of a graph in which any two vertices $X_i$ and $X_j$ are connected by a node if
$S(X_i,R^*_{i,n})\cap S(X_j,R^*_{j,n}) \neq \emptyset$. Since $Q_u = \{1,\ldots,m\}$ is an equivalence class, this graph is connected.
If there is at least one vertex $X_j$ (say) with degree 1, put $L:=j$. Otherwise, the degree of each vertex
is at least two, and we have $m \ge 3$. If $m=3$, the graph is a triangle, and we can choose $L$ arbitrarily. Now suppose the lemma is
true for any graph having $m \ge 3$ vertices, in which each vertex degree is at least 2. If we have an additional $(m+1)$th vertex $X_{m+1}$,
this is connected to at least two other vertices $X_i$ and $X_j$ (say). Of the graph with vertices $X_1,\ldots,X_m$ we can delete one
vertex, and the remaining graph is connected. But $X_{m+1}$ is then connected to either $X_i$ or $X_j$, and we may choose $L=i$ or $L=j$.
Notice that for $d=1$ the proof is trivial since $\cup_{i\in Q_u} S(X_i, R^*_{i,n})$ is an interval, and
we can take either $L=1$ or $L=m$.
\qed

\smallskip
By induction, we now show that
\begin{align}
\label{2*}
\PROB(E_u)
= \textrm{O}\left((\ln n/n)^{|Q_u|-1}\right)
\end{align}
as $n \to \infty$ for each $m:= |Q_u| \in \{2,\ldots, k-1\}$.
We start with the base case $m=2$. Notice that
$
\PROB(E_u)
\le \PROB(\|X_{2}-X_{1}\|\le R^*_{2,n}+R^*_{1,n})
$
and
\begin{align*}
&\PROB(\|X_{2}-X_{1}\|\le R^*_{2,n}+R^*_{1,n}\mid X_1)\\
&=
\PROB(\|X_{2}-X_{1}\|\le R^*_{2,n}+R^*_{1,n},R^*_{2,n}\le R^*_{1,n}\mid X_1)\\
&\quad +
\PROB(\|X_{2}-X_{1}\|\le R^*_{2,n}+R^*_{1,n},R^*_{2,n}>R^*_{1,n}\mid X_1)\\
&\le
\PROB(\|X_{2}-X_{1}\|\le  2R^*_{1,n}\mid X_1)
 +
\PROB(\|X_{2}-X_{1}\|\le 2R^*_{2,n}\mid X_1).
\end{align*}
Now, condition (\ref{L}) entails
\begin{align*}
\PROB(\|X_{2}-X_{1}\|\le  2R^*_{1,n}\mid X_1)
&=
\mu(S(X_1,2R^*_{1,n}))\le
c_{max}\, \mu(S(X_1,R^*_{1,n}))\\
&=
c_{max}\, \frac{y+\ln n}{n}.
\end{align*}
Putting $\widetilde R_{2,n}:=H_{X_2}^{-1}(c_{max}(y+\ln n)/n)$, a second appeal to
(\ref{L}) yields
\begin{align*}
\mu(S(X_2,2R^*_{2,n}))
&\le
c_{max}\, \mu(S(X_2,R^*_{2,n}))
=
c_{max}\, \frac{y+\ln n}{n}
\end{align*}
and thus $2R^*_{2,n}
\le
\widetilde R_{2,n}$.
Consequently,
\begin{align*}
\PROB(\|X_{2}-X_{1}\|\le 2R^*_{2,n}\mid X_1)
&\le
\PROB(\|X_{2}-X_{1}\|\le \widetilde R_{2,n}\mid X_1).
\end{align*}
Let $\gamma_d$ be the minimum number of cones of angle $\pi/3$ centered at $0$ such that their union
covers $\R^d$. Then the cone covering lemma
(cf. Lemma 10.1 in Devroye and Gy\"orfi \cite{DeGy85}, and Lemma 6.2 in Gy\"orfi et al. \cite{GyKoKrWa02}) says that,
 for any $0\le a \le 1$ and any $x_1$, we have
 \begin{align}
\label{cov1}
\mu(\{x_2\in \R^d: \mu(S(x_2,\|x_{2}-x_{1}\|))\le a\})
&\le
\gamma_d \, a.
\end{align}
Now, (\ref{cov1}) implies
\[ 
\mu(\{x_2\in \R^d:\|x_{2}-x_{1}\|\le H_{x_2}^{-1}(a)\})
\gamma_d \, a,
\]
whence
\begin{align*}
\PROB(\|X_{2}-X_{1}\|\le \widetilde R_{2,n}\mid X_1)
&\le
\gamma_d c_{max}\, \frac{y+\ln n}{n}.
\end{align*}
We thus obtain
\begin{align}
\PROB(\|X_{2}-X_{1}\|\le R^*_{2,n}+R^*_{1,n}\mid X_1)
&=\textrm{O}\left( \frac{\ln n}{n} \right),
\label{1*}
\end{align}
and so (\ref{2*}) is proved for $m=2$. For the induction step,
assume (\ref{2*}) holds for $|Q_u|=m\in \{2,\dots , k-2\}$.
If $Q_u$ with $|Q_u|=m+1$ is an equivalence class, then by Lemma \ref{harro} there are random integers $L_1$ and $L_2$ less than $m+2$,
 such that
$Q_u\setminus \{L_1\}$ forms an equivalence class, and
\begin{align*}
\|X_{L_1}-X_{L_2}\|\le R^*_{{L_1},n}+R^*_{L_2,n}.
\end{align*}
It follows that
\begin{align*}
\PROB(E_u)
&\le(m+1)m\PROB\left(E_u\cap\{L_1=m+1,L_2=1\}\right)\\
&\le k(k-1)\PROB\left(\{Q_u\setminus \{m+1\} \mbox{ forms an equivalence class }\}\right. \\
& \left.\quad\cap\{\|X_{m+1}-X_{1}\|\le R^*_{m+1,n}+R^*_{1,n}\}\right)\\
&=
k(k-1)\EXP\big{[}\IND\{Q_u\setminus \{m+1\} \mbox{ forms an equivalence class }\}\\
&\quad\cdot\PROB\left(\|X_{m+1}-X_{1}\|\le R^*_{m+1,n}+R^*_{1,n}\mid X_1,\dots ,X_m\right) \big{]}\\
&=
k(k-1)\EXP\big{[}\IND\{Q_u\setminus \{m+1\} \mbox{ forms an equivalence class }\}\\
&\quad\cdot\PROB\left(\|X_{m+1}-X_{1}\|\le R^*_{m+1,n}+R^*_{1,n}\mid X_1\right)\big{]}\\
&\le
\textrm{O}\left( \frac{\ln n}{n} \right)\PROB\left(Q_u\setminus \{m+1\} \mbox{ forms an equivalence class}\right)\\
&=
\textrm{O}\left( \frac{\ln n}{n} \right)\textrm{O}\left((\ln n/n)^{m-1}\right)\\
&=
\textrm{O}\left((\ln n/n)^{m}\right).
\end{align*}
Notice that the penultimate equation follows from the induction hypothesis, and the last "$\le$" is a consequence of
(\ref{1*}). Notice further that these limit relations imply (\ref{2*}), whence
\begin{align*}
&n^{|Q_u|} \EXP\bigg{[}\ex^{-n\mu\left( \cup_{i\in Q_u} S(X_i, R^*_{i,n})\right)}
\prod_{i\in Q_u}\IND\{\mu (S(X_i,r_{i,k} ))\ge y_n\}
\IND\{E_u\}\bigg{]} \\
&=
\textrm{O}\left(n^{|Q_u|-1-\varepsilon}\right)\PROB(E_u)\\
&=
\textrm{O}\left(n^{|Q_u|-1-\varepsilon}\right) \textrm{O}\left((\ln n/n)^{|Q_u|-1}\right)\\
& = \textrm{o}(1).
\end{align*}
Summarizing, we have shown (\ref{Q}) and thus (\ref{upestim}).
Hence (\ref{cond}) is verified with $\nu =\exp(-y)$, and the theorem is proved.
\qed

\section{Discussion on condition (\ref{INT})}
\label{secdiscuss}

In this final section we comment on condition (\ref{INT}).
For $d=1$, we  verify (\ref{INT}) if on $S(x,r)\cup  S(z,s)$ the distribution function $F$ of $\mu$ is either convex or concave.
 If $\|x-z\| \ge r+s$, then $S(x,r)$ and $ S(z,s)$ are disjoint, therefore suppose
$r+s\ge \|x-z\| \ge \max(r,s)$.
Assume that $F$  is convex, the proof for concave $F$ is similar.
If $x<z$, the convexity of $F$ and
\[
\mu\left( S(z,s)\right)=F(z+ s)-F(z- s)=:p
\]
(say) imply $F(z)-F(z- s)\le p/2$.
Thus
\begin{align*}
\mu\left( S(x,r)\cap  S(z,s)\right)
&=
\mu([z- s,x+ r])\\
&\le
\min\{\mu([z- s,z]),\mu([x,x+ r])\}\\
&=
\min\{F(z)-F(z- s),F(x+ r)-F(x)\}\\
&\le
F(z)-F(z- s)\\
&\le
p/2
\end{align*}
and hence
\begin{align*}
\frac{\mu\left( S(x,r)\cap  S(z,s)\right)}{\mu\left( S(z,s)\right)}
\le \frac 12.
\end{align*}
Thus (\ref{INT}) is satisfied with $\beta =1/2$.

For $d>1$, the problem is more involved.
Again, suppose $r+s\ge \|x-z\| \ge \max(r,s)$.
Writing $\langle \cdot, \cdot \rangle$ for the inner product in $\R^d$, introduce the half spaces
\[ H_1:=\{u \in \R^d: \langle u-x,z -x\rangle \ge 0\}, \quad
H_2:=\{u \in \R^d: \langle u-z,x-z\rangle \ge 0\}.
\]
Then
\begin{align*}
\mu\left( S(x,r)\cap  S(z,s)\right)
&=
\mu\left( (S(z,s)\cap H_2)\cap  (S(x,r)\cap H_1)\right)\\
&\le
\frac{\mu\left( S(z,s)\cap H_2\right)+ \mu\left( S(x,r)\cap H_1\right)}{2}.
\end{align*}
We introduce another implicit condition as follows:
Assume
there are  $\alpha \in (1,2)$ and $\delta >0$ such that, for any $r,s >0$ and any $x,z \in$ supp$(\mu)$ with
$r+s\ge\|x-z\| \ge \max(r,s)$ and $\mu\left( S(x,r)\right)=\mu\left( S(z,s)\right)\le \delta$, one has
either
\begin{align}
\label{a1}
\mu\left( S(z, s)\cap H_2\right)
\le
\alpha \, \mu\left( S(x,r)\cap H_1^c\right)
\end{align}
or
\begin{align}
\label{a2}
\mu\left( S(x,r)\cap H_1\right)
\le
\alpha \mu\left( S(z,s)\cap H_2^c\right).
\end{align}
In case of (\ref{a1}) we have
\begin{align*}
\frac{\mu\left( S(z,s)\cap H_2\right)+ \mu\left( S(x,r)\cap H_1\right)}{2}
&\le
\frac{\alpha \, \mu\left( S(x,r)\cap H_1^c\right)+ \mu\left( S(x,r)\cap H_1\right)}{2}\\
&\le
\alpha\, \frac{ \mu\left( S(x,r)\cap H_1^c\right)+ \mu\left( S(x,r)\cap H_1\right)}{2}\\
&=
\frac{\alpha}{2}\, \mu\left( S(x,r)\right),
\end{align*}
and (\ref{INT}) is verified with $\beta =\alpha/2$. The case of (\ref{a2}) is similar.
For the univariate case and for $x<z$, (\ref{a1}) and (\ref{a2}) mean
\begin{align}
\label{b1}
F(z)-F(z- s)
\le
\alpha \left(F(x)-F(x- r) \right)
\end{align}
and
\begin{align}
\label{b2}
F(x+ r)-F(x)
\le
\alpha \left(F(z+ s)-F(z) \right).
\end{align}
For convex $F$ and small $\delta$, (\ref{b2}) is approximately satisfied with $\alpha \approx 1$. Vice versa, (\ref{b1}) holds  for concave $F$.

\end{document}